\documentclass[11pt]{article}

\usepackage{amssymb,amsthm,amsmath}
\usepackage[utf8]{inputenc}

\usepackage{verbatim}
\usepackage[T1]{fontenc}	
\usepackage[utf8]{inputenc}

\newcommand{\1}{\textbf{1}}

\newcommand{\mb}[1]{\mathbb{#1}}

\newtheorem{thm}{Theorem}
\newtheorem{lem}{Lemma}

\newtheorem{con}{Conjecture}

\theoremstyle{definition}
\newtheorem*{def*}{Definition}
\newtheorem*{question*}{Question}
\newtheorem*{ex*}{Example}

\theoremstyle{remark}
\newtheorem*{rem*}{Remark}

\title{A note on the rational cuspidal curves}

\author{Piotr Nayar, Barbara Pilat}
\date{}

\begin{document}

\maketitle

\begin{abstract}
In this short note we give an elementary combinatorial
argument, showing that the Conjecture of J. Fern\'andez de Bobadilla, I. Luengo, A. Melle-Hern\'andez, A. N\'emethi (see \cite{boro}, Conjecture 1.2) follows from Theorem 5.4 in \cite{boro} in the case of rational cuspidal curves with two critical points. 
\end{abstract}



\section{Introduction}\label{intro}

In this short note we deal with irreducible algebraic curves $C\subset \mathbb{C}P^2$. Such a curve has a finite set of singular points $\{z_i\}_{i=1}^n$ such that a neighbourhood of each singular point intersects $C$ in a cone on a link $K_i \subset S^3$. We would like to know what possible configurations of links $\{K_i\}_{i=1}^n$ arise in this way. We consider only the case in which each $K_i$ is connected (in this case $K_i$ is a knot), and thus $C$ is a rational curve, meaning that there is a rational surjective map $\mb{C}P^1 \to C$. Such a curve is called rational cuspidal. We refer to \cite{moe} for a survey on rational cuspidal curves.

Suppose that $z$ is a cuspidal singular point of a curve $C$ and $B$ is a sufficiently small ball around $z$. Let $\Psi(t)=(x(t),y(t))$ be a local parametrization of $C \cap B$ near $z$. For any polynomial $G(x,y)$ we look at the order at $0$ of an analytic map $t \mapsto G(x(t),y(t)) \in C$. Let $S$ be the set of integers, which can be realized as the order for some $G$. Then $S$ is a semigroup of $\mb{Z}_{\geq 0}$. We call it the semigroup of the singular point, see \cite{wall} for the details and proofs. The gap sequence, $G = \mb{Z}_{\geq 0} \backslash S$, has
precisely $\mu/2$ elements, where the largest one is $\mu-1$. Here $\mu$ stands for the Milnor number. Assume that $K$ is the link of the singular point $z$. The Alexander polynomial of $K$ is of the form
\[
	\Delta_K(t) = \sum_{i=0}^{2m} (-1)^i t^{n_i},
\]
where $(n_i)_{i=0}^{2m}$ form an increasing sequence with $n_0=0$ and $n_{2m}=2g$, for $g$ being the genus of $K$. Writing $t^{2n_i}-t^{2n_{i-1}}=(t-1)(t^{2n_{i}-1}+t^{2n_{i}-2}+\ldots+t^{2n_{i-1}})$ yields the representation
\begin{equation}
\label{deltak}
\Delta_K(t) = 1+(t-1)\sum_{j=1}^k t^{g_j},
\end{equation}
for some finite sequence $0<g_1<g_2<\ldots<g_k$. We have the following lemma (see \cite{wall}, Exercise 5.7.7), which relates the Alexander polynomial to the gap sequence of a singular point. 

\begin{lem} The sequence $g_1,\ldots, g_k$  in (\ref{deltak}) is the gap sequence of the semigroup of the
singular point. In particular, $k = |G| = \mu/2$, where $\mu$ is the Milnor number, so $|G|$
is the genus.
\end{lem}
If we write $t^{g_j}= (t-1)(t^{g_j-1}+t^{g_j-2}+\ldots+1)+1$, we obtain
\[
	\Delta_K(t) = 1+(t-1)g(K) + (t-1)^2 \sum_{j=0}^{\mu-2} k_j t^j,
\]
where $k_j = |\{ m>j: m \notin S \}|$. This motivates the following definition. 

\begin{def*} For any finite increasing sequence of positive integers $G$ we define
\[
	I_G(m) = |\{k \in G \cup \mb{Z}_{<0}: \  k \geq m\}|,
\]
where $\mb{Z}_{<0}$ is the set of negative integers. We shall call $I_G$ the gap function, because in most applications $G$ will be a gap sequence of some semigroup.
\end{def*}

\noindent Clearly, for $j=0,1,\ldots, \mu-2$ we have $I_G(j+1)=k_j$. 

In \cite{1} the following conjecture was proposed.

\begin{con}\label{conjecture}
Suppose that the rational cuspidal curve $C$ of degree $d$ has critical points $z_1, \ldots, z_n$. Let $K_1, \ldots, K_n$ be the corresponding links of singular points and let $\Delta_1, \ldots, \Delta_n$ be their Alexander polynomials. Let $g$ be the genus of $K$. Let $\Delta = \Delta_1, \ldots, \Delta_n$, expanded as 
\[
	\Delta(t) = 1 + \frac{(d-1)(d-2)}{2}(t-1) + (t-1)^2 \sum_{j=0}^{2g-2} k_l 
\]
Then for any $j=0,\ldots, d-3 $ we have $k_{d(d-j-3)} \leq (j+1)(j+2)/2$, with equality for $n=1$.
\end{con}
\noindent This conjecture was verified in the case $n=1$ by Borodzik and Livingston, see \cite{boro}.

We define the infimum convolution of two functions.

\begin{def*}
Let $I_1,I_2,\ldots,I_n: \mb{Z} \to \mb{Z}_{\geq 0}$. We define
\[
	(I_1 \diamond I_2 \ldots \diamond I_n)(k) = \min_{\begin{array}{cc}
	\scriptstyle{k_1, k_2,\ldots,k_n \in \mb{Z}} \\
	\scriptstyle{k_1+k_2+\ldots+k_n = k}
	\end{array}}\left( I_1(k_1) + I_2(k_2)+\ldots+  I_n(k_n)  \right).
\]	 
\end{def*} 

\noindent In \cite{boro} the authors gave the proof of the following theorem.

\begin{thm} (\cite{boro}, Theorem 5.4) Let $C$ be a rational cuspidal curve of degree $d$. Let $I_1,\ldots,I_n$ be the gap functions associated to each singular point on $C$. Then for any $j \in \{-1,0,\ldots, d-2\}$ we have
\[
	I_1 \diamond I_2 \diamond \ldots \diamond I_n(jd+1) = \frac12 (j - d + 1)(j - d + 2).
\]
\end{thm}
   
\noindent Note that we have $|G_1|+ |G_2|+\ldots+|G_n|=\frac{(d-1)(d-2)}{2}$. Therefore, one can give an equivalent reformulation of the Conjecture \ref{conjecture}. 

\begin{con}\label{conjecture2}
Suppose that the rational cuspidal curve $C$ of degree $d$ has critical points $z_1, \ldots, z_n$. Let $K_1, \ldots, K_n$ be the corresponding links of singular points and let $\Delta_1, \ldots, \Delta_n$ be their Alexander polynomials. Moreover, let $G_1, G_2,\ldots,G_n$ be the gap sequences of these points. Let $g=|G_1|+ |G_2|+\ldots+|G_n|$ be the genus of $K$. Let $\Delta = \Delta_1, \ldots, \Delta_n$, expanded as 
\[
	\Delta(t) = 1 + (t-1)g + (t-1)^2 \sum_{j=0}^{2g-2} k_l 
\]
and let $I = I_1 \diamond I_2 \diamond \ldots \diamond I_n$.
Then for any $j=0,\ldots, d-3$ we have $k_{d(d-j-3)} \leq I(d(d-j-3)+1)$, with equality for $n=1$.
\end{con}

In this note we give an elementary argument, showing that FLMN conjecture follows from \cite{boro} for $n=2$. The idea of our proof is to forget about the specific structure of the problem coming from theory of singularities and to prove Conjecture \ref{conjecture2} for general sets $G_1, G_2$. Namely, we have the following theorem. 

\begin{thm}
Let $G, H$ be two finite sets of positive integers and let $I_G,I_H:\mb{Z} \to \mb{Z}_{\geq 0}$ be their gap functions. Let us define the polynomials
\[
	\begin{array}{l}
		\Delta_G(t)   =  1 + (t-1) \sum_{j=1}^{|G|} t^{g_j}   =  1+(t-1)|G| + (t-1)^2 \sum_{j \geq 0} k_j^G t^j \\
		\Delta_H(t)   = 1 + (t-1) \sum_{j=1}^{|H|} t^{h_j}  = 1+(t-1)|H| + (t-1)^2 \sum_{j \geq 0} k_j^H t^j
	\end{array},
\] 
where $k_j^G = I_G(j+1), k_j^H = I_H(j+1)$, $j \geq 0$. Take $\Delta = \Delta_G \cdot \Delta_H$ and $I =  I_G \diamond I_H$. Then 
\[
	\Delta(t) = 1+ (t-1)(|G|+|H|) + (t-1)^2 \sum_{j \geq 0} k_j t^j,
\]
where $k_j \leq I(j+1)$ for $j \geq 0$.
\end{thm}

\noindent This gives the proof of Conjecture \ref{conjecture} in the case $n=2$.

It is natural to ask whether the above theorem is valid for arbitrary $n \geq 2$. Recently, after we found our elementary combinatorial argument for $n=2$, J. Bodn\'ar and A. N\'emethi showed that the Conjecture \ref{conjecture} is false for $n \geq 3$, see \cite{bn}. They also found yet another proof of Conjecture \ref{conjecture} in the case of two singularities.

\section{Proof of the main result}\label{proof}

In this section we give a proof of our main result.

\begin{proof}
Our goal is to express the numbers $k_j$ in terms of the numbers $k_j^G$ and $k_j^H$. We have
\begin{align*}
& \Delta(t)  = \Delta_G(t) \Delta_H(t) = 1+ (t-1)(|G|+|H|) \\ &  + (t-1)^2 \Big[ |G|\cdot |H| + \sum_{j \geq 0} (k_j^G + k_j^H)t^j + (t-1)(|G| \sum_{j \geq 0} k_j^H t^j + |H| \sum_{j \geq 0} k_j^G t^j  ) \\ & + (t-1)^2 \Big( \sum_{j \geq 0} k_j^G t^j  \Big) \Big( \sum_{j \geq 0} k_j^H t^j  \Big) \Big] =  1+ (t-1)(|G|+|H|) + (t-1)^2 \Theta(t), 
\end{align*}
with
\[
	\Theta(t)  = |G|\cdot |H| +k_0^G(1-|H|) + k_0^H(1-|G|) + k_0^G k_0^H  
	 +  \sum_{j \geq 1} t^j  k_j,  
\]
where
\[
	k_j =  k_j^G(1-|H|) +|H| k_{j-1}^G + k_j^H(1-|G|) +|G| k_{j-1}^H +l_j	  
\]
and
\[
	 l_j = \sum_{   u+v=j,\  u,v \geq 0 } k_u^G k_v^H - 2\sum_{u+v=j-1,\  u,v \geq 0} k_u^G k_v^H + \sum_{u+v=j-2,\  u,v \geq 0} k_u^G k_v^H.
\]
Note that $k_0^G = |G|$ and $k_0^H = |H|$.  Therefore,
\begin{align*}
	k_0 & = |G|\cdot |H| +k_0^G(1-|H|) + k_0^H(1-|G|) + k_0^G k_0^H \\ & = |G|\cdot |H| + |G|(1-|H|) +|H|(1-|G|) + |G| \cdot |H| = |G|+|H|.
\end{align*}
Moreover, for $k \geq 1$ we have
\[
	I_H(k) \geq |H|-(k-1), \qquad I_G(1-k) = |G|+(k-1). 
\]
Thus, for $k \geq 1$ we obtain
\[
	I_G(1-k) + I_H(k) \geq |G|+|H|.
\]
For $k \leq -1$ 
\[
	I_G(1-k) \geq |G|+k, \qquad I_H(k) = |H|-k.
\]
In this case we obtain
\[
	I_G(1-k) + I_H(k) \geq |G|+|H|
\]
and we arrive at
\begin{align*}
I(1) & = \min_{k \in \mb{Z}}(I_G(1-k)+I_H(k)) \\ & = I_G(1)+I_H(0) = I_G(0)+I_H(1) = |G|+|H|= k_0. 
\end{align*}

Note that
\begin{align*}
	 l_j & = \sum_{   u+v=j,\  u,v \geq 0 } k_u^G k_v^H - \sum_{u+v=j,\  u \geq 0,v \geq 1} k_u^G k_{v-1}^H - \sum_{u+v=j,\  u \geq 1,v \geq 0} k_{u-1}^G k_{v}^H \\
	 & \quad + \sum_{u+v=j,\  u,v \geq 1} k_{u-1}^G k_{v-1}^H = \sum_{u+v=j, u,v \geq 1} (k_u^G-k_{u-1}^G)(k_v^H-k_{v-1}^H) \\
	 & \quad  + k_0^G k_j^H + k_j^G k_0^H - k_0^G k_{j-1}^H -k_{j-1}^G k_0^H .
\end{align*}
Thus,
\begin{align*}
	k_j & = \sum_{u+v=j, u,v \geq 1} (k_u^G-k_{u-1}^G)(k_v^H-k_{v-1}^H) \\
	& \qquad + \Big( (k_0^G k_j^H + k_j^G k_0^H - k_0^G k_{j-1}^H -k_{j-1}^G k_0^H) \\
	& \qquad +( k_j^G(1-|H|) +|H| k_{j-1}^G + k_j^H(1-|G|) +|G| k_{j-1}^H )  \Big).
\end{align*}
Observe that we have a miracle,
\begin{align*}
	 k_0^G k_j^H & + k_j^G k_0^H - k_0^G k_{j-1}^H -k_{j-1}^G k_0^H \\ & \hspace{2cm} + k_j^G(1-|H|) +|H| k_{j-1}^G + k_j^H(1-|G|) +|G| k_{j-1}^H  \\
	&  = |G| k_j^H + k_j^G |H| - |G| k_{j-1}^H -k_{j-1}^G |H| \\
	&  \hspace{1cm} + k_j^G(1-|H|) +|H| k_{j-1}^G + k_j^H(1-|G|) +|G| k_{j-1}^H \\
	&  = k_j^G + k_j^H.
\end{align*}
We get
\[
	k_j = k_j^G + k_j^H + \sum_{u+v=j, u,v \geq 1} (k_{u-1}^G-k_{u}^G)(k_{v-1}^H-k_{v}^H).
\]

We are to prove that $k_j \leq (I_G \diamond I_H)(j+1)$. It suffices to prove that $k_j \leq I_G(j+1-l)+I_H(l)$ for every $l \in \mb{Z}$. Thus, we have to deal with the inequality
\[
	 k_j^G + k_j^H + \sum_{u+v=j, u,v \geq 1} (k_u^G-k_{u-1}^G)(k_v^H-k_{v-1}^H) \leq  I_G(j+1-l)+I_H(l), \quad j \geq 1, l \in \mb{Z}. 
\]
Note that if $u+v=j$ then we have either $u \geq j-l+1$ or $v \geq l$. Thus,
\[
	\1_{u \in G} \1_{v \in H} \1_{u+v=j} \leq \1_{u \in G \cap [j-l+1,j]} + \1_{v \in H \cap[l,j]}. 
\]
In the above expression we have used the convention $[a,b]=\emptyset$ for $a>b$.
We obtain
\begin{align*}
	\sum_{u+v=j, u,v \geq 1} & (k_{u-1}^G-k_{u}^G)(k_{v-1}^H-k_{v}^H)  = \sum_{u+v=j, u,v \geq 0} (k_{u-1}^G-k_{u}^G)(k_{v-1}^H-k_{v}^H) \\
	& = \sum_{u+v=j, u,v \geq 0} \1_{u \in G} \1_{v \in H} \leq \sum_{u+v=j, u,v \geq 0} \left( \1_{u \in G\cap[j-l+1,j]} + \1_{v \in H\cap[l,j]} \right) \\ & = (k_{j-l}^G-k_j^G) + (k_{l-1}^H-k_j^H), 
\end{align*}
what finishes our proof.
\end{proof}

\section*{Acknowledgements}

We are grateful to Maciej Borodzik for his valuable comments about the meaning of our result.

\vspace{1cm}

\noindent Piotr Nayar, \texttt{nayar@mimuw.edu.pl} \\
\noindent Institute of Mathematics, University of Warsaw, \\
\noindent Banacha 2, \\
\noindent 02-097 Warszawa, \\
Poland. 

\vspace{1em}

\noindent Barbara Pilat, \texttt{B.Pilat@mini.pw.edu.pl} \\
\noindent Faculty of Mathematics and Information Science, Warsaw University of Technology, \\
\noindent Koszykowa 75, \\
\noindent 00-662 Warszawa, \\
Poland.


\begin{thebibliography}{9}

\bibitem[FLMN]{1} J. Fern\'andez de Bobadilla, I. Luengo, A. Melle-Hern\'andez, A. N\'emethi, \emph{Classification of rational unicuspidal projective curves whose singularities have one Puiseux pair}, Proceedings of Sao Carlos Workshop 2004 Real and Complex Singularities, Series Trends in Mathematics, Birkh\"auser 2007, 31–46.

\bibitem[BL]{boro} M. Borodzik, C. Livingston, \emph{Heegaard Floer homology and rational cuspidal curves}, arXiv:1304.1062 

\bibitem[BN]{bn} J. Bodn\'ar and A. N\'emethi, arxiv:1405.0437

\bibitem[M]{moe} K. Moe, \emph{Rational cuspidal curves}, Master Thesis, University of Oslo 2008, permanent link at University of Oslo: \\ https://www.duo.uio.no/handle/123456789/10759

\bibitem[W]{wall} C. T. C. Wall, \emph{Singular Points of Plane Curves}, London Mathematical Society Student Texts, 63. Cambridge University Press, Cambridge, 2004.

\end{thebibliography}
\end{document}